\documentclass[10pt,twoside]{article}
\usepackage{graphicx}
\usepackage{amsmath, amssymb}
\usepackage{Latex-document}

\newtheorem{theorem}{Theorem}
\newtheorem{proposition}[theorem]{Proposition}
\newtheorem{lemma}[theorem]{Lemma}

\newtheorem{fact}[theorem]{Fact}

\newcommand{\qed}{\relax\ifmmode\eqno\B ox\else\mbox{}\quad\nolinebreak\hfill$\Box$\smallskip\fi}

\newcommand{\lclosed}{$\lambda$-closed }
\newcommand{\Lpinf}{L^{p^{\infty}}}

\newcommand{\Apinf}{p^{\infty} A(L)}

\newcommand{\infdef}{ $\infty$-definable }

\newcommand{\K}{K[X_{\infty}]}
\newcommand{\KX}{\K}

\newcommand{\pupu}{(p^\nu)^n}

\newcommand{\pomega}{(p^\nu)^{< \omega}}

\title{\bf  Groups Interpretable   \vskip -2mm
in Theories of Fields \vskip 6mm}

\author{E. Bouscaren\vspace*{-0.5cm}\thanks{University Paris 7 -
    CNRS, Department of Mathematics, Case 7012, 2 Place Jussieu,
    75251 Paris Cedex 05, France. E-mail: elibou@logique.jussieu.fr}}
\date{\vspace{-8mm}}

\begin{document}

\maketitle

\thispagestyle{first} \setcounter{page}{3}

\begin{abstract}\vskip 3mm  We  survey some   results  on the
  structure of the groups which are definable in
theories of fields involved in  the applications of
model theory to Diophantine geometry. We focus more particularly on
separably closed fields of finite degree of imperfection.
\vskip 4.5mm


\end{abstract}

\vskip 12mm

\section{Introduction} \label{section1}\setzero
\vskip-5mm \hspace{5mm}

In the last ten years, the model theory of fields has seen
striking new developments,  with applications in particular
to differential algebra and Diophantine geometry. One of the main
ingredients in these applications is the analysis of the structure
of groups definable  in fields with  added ``definable structure''.

Model theory studies structures  with a family of
distinguished subsets of  their Cartesian products,
the family of {\em definable} subsets, which is requested to be
closed under finite Boolean operations and projections.
In the case of algebraically closed fields,
the  definable sets are exactly the constructible sets
in the Zariski topology (finite Boolean combinations of Zariski
closed sets).
If one considers fields which are not algebraically closed (for
example,  fields of positive characteristic which are
separably closed and not
perfect) or algebraically closed fields with new
operators (differentially closed fields, fields with a generic
automorphism), then the family of definable sets is much richer
than the family of Zariski constructible sets. In each of the above
cases, one
can generalize the classical geometric notions, by using the
tools  developed by model theory (abstract notion of independence,
of dimensions...). For example:

1.  One can define ``good'' topologies
which strictly contain the Zariski topology.

2. Different  notions of dimensions can be attached to
definable sets (or  infinite
intersections of definable sets, which we call {\em infinitely
  definable, or $\infty$-definable,  sets}). In the case of algebraically
closed fields, all such  possible notions of abstract dimension must
coincide
and be equal to the classical algebraic dimension. In the other cases,
these dimensions may be different, some may take infinite ordinal
values or may be defined only
for some special classes of definable (or $\infty$-definable)
sets.

3. If $K$ is any of the above mentioned fields, and if $H$ is an
algebraic group defined over $K$, then the group $H(K)$ of the
$K$-rational points of $H$ is a definable group. But there are
``new'' families of definable groups which are not of this form.

In fact, it is precisely the study of  certain specific families of ``new''
definable groups of finite dimension which are at the center
of the applications to
Diophantine geometry. We will not attempt here to explain how
the model theoretic
analysis of the finite rank definable groups yields these
applications.  There have been in recent years many surveys
and presentations of the subject to which we refer the reader (see
for example, \cite{B1},\cite{B2}, \cite{Hr0}, \cite{Pi3} or
\cite{scanlon4}). We will come back to this subject, but  very
briefly,  at the
end  in Section \ref{section4}.

The first general question raised by the existence of these new
definable groups is that of
their relation to the classical algebraic groups. Remark that this  question
already makes sense in the context of ``pure'' algebraically closed
fields, about the class of definable (= constructible) groups. In that case,
it is true that any constructible group in an algebraically
closed field $K$ is  constructibly isomorphic to
the $K$-rational points of an algebraic group defined over $K$ (see
for example \cite{bousweil} or \cite{pillaymanch}).

Let us now consider briefly the case of a field $K$
of characteristic $p>0$ which is separably closed and not perfect.
Then the class of constructible sets is no longer closed under
projection and there  are many  definable groups which are not
constructible, the most obvious one  being  $K^{p}$. There are also some
groups which are proper intersections of infinite  descending chains of
definable  groups:  for example, $ K^{p^\infty}(=\bigcap_n K^{p^n})$,
the field of infinitely
$p$-divisible elements of the multiplicative group,
or $\bigcap_n p^n A(K)$, for $A$ an
Abelian variety defined over $K$.

It is nevertheless true, as we will see,   that
every definable group in $K$ is definably isomorphic to the $K$-rational
points of an algebraic group defined over $K$. Furthermore, as
in the classical case of one-dimensional  algebraic groups, it is
possible to give a complete description, up to definable isomorphism,
of the one-dimensional infinitely definable
groups.

There are results of similar type for the other classes of
enriched fields mentioned above. In  this short paper, we will
concentrate mainly on the case of
separably closed fields (in Section \ref{section3}). Before this,
in Section \ref{section2}, we will only   very briefly
present the model theoretic  setting for two other examples of
``enriched'' fields,  in
characteristic zero, differentially closed fields and generic
difference fields. We hope this will give  the reader an idea of
what the common
features and the differences might be in the model theoretic
analysis of these different classes of fields.

Finally, there are of course many other classes of fields whose model
theory has been extensively developed in the past years with many
connections to algebra, semi-algebraic or subanalytic geometry, and which we
are not going to mention here: for example, valued fields, ordered
fields,
``o-minimal'' expansions of the real field...

\section{Two short examples}\label{section2}\setzero
\vskip-5mm \hspace{5mm}

We will just very briefly describe the two characteristic zero
examples mentioned above.

\subsection{Differentially closed fields of characteristic zero}\label{DCF}

\vskip-5mm \hspace{5mm}

We consider a field $K$ of characteristic zero, with a {\em
  derivation} $\delta$, that is,   an additive map from $K$ to $K$
  which satisfies that for all $x,y$ in $K$, $\delta(xy) =
x\delta(y)+y \delta(y)$. We define the ring $K_\delta [X]$ of {\em differential polynomials}  over $K$ to be the
ring of polynomials in infinitely many variables $K[X,\delta(X),\delta^2(X),\cdots,\delta^n (X),\cdots]$. The {\em
order} of the differential polynomial $f(X)$ in $K_\delta[X]$ is $-1$ if $f \in K$ and otherwise the largest $n$
such that $\delta^n(X)$ occurs in $f(X)$ with non zero coefficient. We say that $K$ is {\em differentially closed}
if for any non-constant differential polynomials $f(X)$ and $g(X)$, where the order of $g$ is strictly less than
the order of $f$, there is a $z$ such that $f(z) = 0$ and $g(z) \not= 0$. In model theoretic terms, this means
exactly
  that $K$ is existentially closed.

From now on we suppose that $(K,\delta)$ is a large differentially
closed field (a universal domain).

We say that $F\subseteq K^n$ is a {\em  $\delta$-closed set},  if there are $f_1,\cdots, f_r \in K_\delta
[X_1,\cdots,X_n]$ such that $F = \{ (a_1,\cdots,a_n) \in K^n ; f_1(a_1,\cdots,a_n) =\cdots = f_r(a_1,\cdots,a_n) =
0\}$. The ring $K_\delta [X_1,\cdots,X_n] $ is of course not Noetherian but the $\delta$-closed sets (which
correspond to radical differential ideals) form the closed sets of a Noetherian topology on $K$, the {\em
$\delta$-topology}.

We now consider the {\em $\delta$-constructible} sets, that is, the finite
Boolean combinations of $\delta$-closed sets. This class is closed under
projection (this is quantifier elimination for the theory),
hence the definable sets (we call them
{\em $\delta$-definable} sets) are
exactly the $\delta$-constructible sets.
To every $\delta$-definable set one can
associate a dimension (the Morley rank) which can take infinite
countable ordinal values.

There are ``new'' definable groups, which are not of the form
$H(K)$ for any  algebraic group $H$. In particular, any $H(K)$
will have infinite dimension. In contrast, the
{\em  field of constants of $K$}, ${ Cons(K)}= \{a\in K; \delta(x)
= 0\}$, is a $\delta$-closed
set which is not constructible;  it is  an algebraically closed
subfield of $K$ and has dimension  one.

Nevertheless the following is true:
\begin{proposition}\label{diff} {\em (\cite{pillay0})}
Let $G$ be a  $\delta$-definable group in $K$.  Then there is an
algebraic group $H$, defined over $K$,  such that $G$ is
definably isomorphic to
a ($\delta$-definable) subgroup of $H(L)$.
\end{proposition}

For the many more existing results on  $\delta$-definable groups,
we refer the reader to \cite{pillaydiff}, or from the differential
algebra point of view, to \cite{buium}.

\subsection{Generic difference fields}\label{ACFA}

\vskip-5mm \hspace{5mm}

We now consider an algebraically closed field $K$
 with an automorphism $\sigma$. We say that $(K,\sigma)$ is
a generic difference field if every difference equation which has a
solution in an extension of $K$ has a solution in $K$. The theory
of generic difference fields  has been extensively studied in
\cite{zoeudi} and \cite{zoeudikobi}.

Let us suppose that $(K,\sigma)$ is a generic difference field in characteristic zero. We consider the ring of
{\em $\sigma$-polynomials},
$$ K_\sigma [X_1,\cdots,X_n] = K[X_1,\cdots,X_n,\sigma(X_1),: \cdots,\sigma(X_n),\sigma^2(X_1),\cdots,
\sigma^2(X_n),\cdots]. $$
 We say that $F\subseteq K^n$ is a
{\em f $\sigma$-closed set} if there are $f_1,\cdots, f_r \in K_\sigma [X_1,\cdots,X_n]$ such that $F = \{
(a_1,\cdots,a_n) \in K^n :  f_1(a_1,\cdots,a_n) = \cdots = f_r(a_1,\cdots,a_n) = 0\}$.  The $\sigma$-closed sets
form the closed sets of a Noetherian topology on $K$, the {\em $\sigma$-topology}. The class of $\sigma$-definable
sets is the closure under finite Boolean operations and
 projections
of the $\sigma$-closed sets.

Again there are ``new'' $\sigma$-definable groups. For example, the field $Fix(K) = \{a\in K : \sigma(a) = a\}$,
the fixed field of $\sigma$ in $K$, is a $\sigma$-closed set of dimension one.

Here the best result possible for arbitrary $\sigma$-definable groups is the following:

\begin{proposition}\label{gen}{\em (\cite{piko})} Let $G$ be a group definable in
  $(K,\sigma)$. Then there are an algebraic group $H$ defined over
  $K$, a finite normal subgroup $N_1$ of $G$,
  a $\sigma$-definable subgroup $H_1$ of $H(K)$  and a finite
  normal subgroup $N_2$ of $H_1$, such that $G/N_1$ and $H_1/N_2$
  are $\sigma$-definably isomorphic.\end{proposition}

The analysis of groups of finite dimension is one of the main tools
in Hrushovski's proof of the Manin-Mumford conjecture in
\cite{Hr2}.

\section{Separably closed fields of finite degree of imperfection}\label{section3}
\setzero\vskip-5mm \hspace{5mm}

Separably closed fields are particularly interesting from the model
theoretic point of view for many reasons, in addition to the fact that
they form  the framework for Hrushovski's proof of the
 Mordell-Lang conjecture in charactersitic $p$.  Let us just mention one
reason
here: they are the only fields known to be  stable and non
superstable, and in fact it is conjectured that they are the only
existing ones.

We will just focus on the main properties of the groups that are
definable in a separably closed field of finite degree of
imperfection, but we need first to introduce some notation and recall
some basic facts (see \cite{delon2}).

\subsection{Some basic facts and notation}

\vskip-5mm \hspace{5mm}

Let $L$ be a separably closed field of charcteristic $p>0$ and
of finite degree of
imperfection which is not perfect, i.e.,  $L$ has no proper separable
algebraic extension, and  $|L:L^p| =
p^\nu$, with $0 <\nu$. In order to avoid confusion we denote the
Cartesian product of $k$ copies of $L$ by $L^{\times k}$.

A subset $B = \{b_1,\cdots,b_\nu\}$ of $L$ is called a {\em
  $p$-basis}  of $L$ if the {\em set of $p$-monomials of $B$},
$\{M_j := \prod_{i=1}^{\nu}{b_i}^{j(i)} ; j \in p^\nu \}$
forms a linear basis of $L$ over $L^p$.
Each element $x$ in $L$ can be written in a unique way as
$x = \sum_{j \in p^\nu} {x_j}^{p} M_j$.
{\bf From now on  we fix a $p$-basis $B$ of $L$ and the $M_j$'s, with $j\in
p^\nu$, always denote the $p$-monomials of $B$}. We  suppose that
$L$ is large (a universal domain, or in model theoretic terms, saturated)
and we fix some small separably closed subfield $K$ of $L$,
containing $B$ and of same degree of imperfection $\nu$.

We let  $f_j$ denote the map  which to $x$ associates $x_j$.  The $x_j$'s are called the {\em $p$-components of
$x$ of level one}. More generally, one can associate to $x$ a tree of countable height indexed by
$(p^{\nu})^{<\omega}$, which we call the {\em tree of $p$-components of $x$}. For $\sigma \in
(p^{\nu})^{<\omega}$, we define $x_\sigma$ by induction: $x_\emptyset = x$  and if $\tau \in(p^{\nu})^n$, and $j
\in p^\nu$, we let $x_{(\tau,j)}$ be equal to $f_j(x_\tau )$; $x_{(\tau,j)}$ is called a $p$-component of $x$ of
level $n+1$.

We will also use the notation $a_\infty := {(a_\sigma)}_{\sigma \in (p^{\nu})^{<\omega}}$, for $a \in L$.

\noindent{\bf The ring \boldmath$\K$.} $\K$ is the polynomial ring in countably many indeterminates indexed in a
way which will allow the natural substitution by the $p$-components of elements: for $X$ a single variable,
$X_\infty:=(X_\sigma)_{\sigma \in \pomega}$, and for $X=(Y_1,\dots,Y_k)$ a $k$-tuple of variables, $X_\infty:=
((Y_1)_\infty,\dots,(Y_k)_\infty)$. The  ring $K[X_{\infty}]$ is a countable union of Noetherian rings, hence each
ideal is countably generated. We let $I^0(X)$ denote the ideal of $\KX$ generated by the polynomials $X_\sigma
-\sum_{j \in p^{\nu}} X_{(\sigma,j)}^p M_j$, $\sigma \in (p^\nu)^{< \omega}$.

\subsection{The \boldmath$\lambda$-topology}\label{language}

\vskip-5mm \hspace{5mm}

Given a set of polynomials $S$ of $\K$, let $ V(S)=\{a\in L^{\times k} : f(a_\infty) = 0 \hbox{ for all } f\in S
\}$. Such a $V(S)$ is called $\lambda$-{\em {closed}} (with parameters in $K$ or over $K$) in $L$.

Given $A\subseteq L^{\times k}$, we define its {\em canonical ideal } $I(A)$ over $K$, $I(A):=\{f\in \K :
f(a_\infty)=0 \hbox{ for all } a \in A\}$.

The $\lambda$-closed subsets of $L^{\times k}$ form  the closed sets of the {\em $\lambda$-topology} on $L^{\times
k}$. This topology is not Noetherian but is the limit of countably many Noetherian topologies.

Let $C$ be a commutative $K$-algebra.   An ideal $I$ of $C$
is  {\em separable}
if, for all $c_j \in C$, $j \in p^{\nu}$,
$\hbox{if }\sum_{j \in p^{\nu}} c_j^p M_j \in I ,
 \hbox{ then each } c_j \in I$.

\begin{fact}[``Nullstellensatz'']\label{delon-nullstell}
1. The map $A\mapsto I(A)$ induces a bijection between $\lambda$-closed subsets of the affine space $L^{k}$ which
are defined over $K$ and ideals of $\K$ which are separable and contain $ I^0(X)$. The inverse map is $I\mapsto
V(I)$.
\end{fact}

Now for the basic properties of the first-order theory:

\begin{fact}\label{basicfacts} 1. The theory of separably closed fields of
characteristic $p$, of degree of imperfection $\nu$,
and with $p$-basis $\{ b_1,...,b_{\nu} \}$ is complete and
admits elimination of quantifiers and
elimination of imaginaries in the language
$$ {\cal L} _{p,\nu} = \{ 0,1,+,-,. \}
\cup \{ b_1,...,b_{\nu} \} \cup \{ f_i ; i \in p ^{\nu} \}.$$
\end{fact}

In particular, any definable set is $\lambda$-constructible, that is, a finite Boolean combinations of definable
$\lambda$-closed sets. Remark that it is impossible to associate to an arbitrary definable set a well-behaved
notion of dimension: indeed, such a dimension would need to be invariant under definable bijections, but for every
$n$ the map $\lambda_n$, defined by $\lambda_n(x) := (x_\sigma )_{\sigma  \in \pupu}$, is a definable bijection
between $L$ and $L^{\times p^{\nu n}}$. But some \infdef  sets will have a well-defined dimension, for example the
field $\Lpinf := \bigcap_n L^{p^n}$,  which is the biggest algebraically closed subfield of $L$, has dimension
one. In fact, $\Lpinf$ is the unique (up to definable isomorphism) infinitely definable field of dimension one
(\cite{margit}, \cite{Hr1}).

\subsection{Definable groups}

\vskip-5mm \hspace{5mm}

Again, amongst the definable groups, one finds the ``classical''
ones, that is groups of the form $H(L)$ for $H$ any algebraic group
defined over $L$. These groups have certain specific properties which are
not true of all the definable groups in $L$. Recall
that a definable subset
$X$ of $G$ is said to be {\em generic}
if $G$ is covered by a finite number of translates of $X$, and an
element of $G$ is {\em generic for the group} if every definable
set which contains it is generic. In an algebraic group, generics in
the topological sense coincide with generics for the algebraic group. Recall
also that a definable group is said to be {\em connected} if it has no
proper definable subgroup of finite index, and {\em
  connected-by-finite} if it has a definable connected subgroup of
finite index.

\vspace{-0.2cm}
\begin{proposition}\label{alggroup}{\em (\cite{BD1}, \cite{Hr1})} Let $H$ be an algebraic group defined over
$K$. Then $H(L)$ is connected-by-finite. If  $H$ is connected
(hence irreducible as an
algebraic group), then $H(L)$ is connected (and irreducible for the
$\lambda$-topology) and if $a\in H(L)$ is a
generic point, then the ideal $I(a) = \{f\in \K : f(a_\infty) =0\}$
is minimal amongst the ideals $I(h)$, for $h \in
H(L)$. \end{proposition}

The above says that in the group $H(L)$, the generics in the topological sense coincide with generics for the
group. In an arbitrary group defined in $L$, this need not be the case.

Consider the definable bijection $f$ from $L$ to $L$ defined in the following way: if $x \in L\setminus L^p$,
$f(x) = x^p$, if $x \in L^p\setminus L^{p^2}$, $f(x) = x^{1/p}$, if $x \in {L^{p^2}}$, $f(x) = x$.

Transporting addition through $f$, one gets a group on $L$ again, $G := (L,*)$, definably isomorphic to $(L,+)$,
hence connected. The set $L$ itself is of course $\lambda$-closed and irreducible with associated ideal $I(L) =
I^0(X)$. The ideal associated to the (group) generic of $(L,*)$ is generated by  $I^0(X)$ and $\{X_i=0 : i \in
p^\nu , i \not= 0\}$, and strictly contains $I^0(X)$.

This question of the uniqueness of the notion of generic is not the only one posing problems for arbitrary
definable groups in $L$. For example, there is no reason, coming  from general properties of stable (non
superstable) theories,  which a priori forces all these definable groups to be connected-by-finite.

Nevertheless, one can in fact show that the situation is as close to the classical one as it could be:
\begin{proposition}\label{defgroup}{\em \cite{BD1}} Every definable
  group $G$ in $L$  is connected-by-finite and is  definably
  isomorphic to the group of $L$-rational points of an  algebraic
  group $H$ defined over $L$.
\end{proposition}

One more remark, in the case of algebraic groups, by Prop. \ref{alggroup},  irreducibility transfers down to the
set of $L$-rational points. But this is not the case for an arbitrary variety: if one considers for example the
irreducible variety defined by the equation $Y^{p^m}X + Z^{p^m} = 0$, for $m \geq 1$, then the $\lambda$-closed
set $V(L)$ is no longer  irreducible in the sense of the $\lambda$-topology.

\subsection{Minimal groups}

\vskip-5mm \hspace{5mm}

The previous result enables us to give a complete description of groups of dimension one, and more generally of
some classes of commutative groups.

We say that an $\infty$-definable set  $D$ is {\em minimal} if any definable subset of $G$ is finite or co-finite.
If $D$ is actually definable, then we say that $D$ is {\em strongly minimal}.

The minimal groups are exactly the connected groups of dimension
 (U-rank) equal to one. A minimal group must be commutative.

From the basic properties of commutative algebraic groups over an algebraically closed field of characteristic $p$
and Proposition \ref{defgroup}, one can deduce:
\begin{lemma} Let $G$ be a minimal group $\infty$-definable in
  $L$, then $G$ has exponent $p$ or $G$ is divisible. \end{lemma}

We first consider the commutative groups of exponent $p$:

\begin{proposition} {\em \cite{BD2}} Let $G$ be a commutative \infdef group of
  exponent $p$ definable in  $L$. Then $G$ is definably isomorphic
  to a $\lambda$-closed subgroup of the additive
  group $(L,+)$. Furthermore, if $G$ is definable, then it is definably isogenous to the
  group of $L$-rational points of a vector group. \end{proposition}

Note that even when $G$ is connected it is not necessarily definably isomorphic to the group of rational points of
a vector group.

Then we consider the commutative divisible groups, which we show to be exactly the ones that were considered by
Hrushovski in  \cite{Hr1}:

\begin{proposition} {\em \cite{BD2}} 1. Let $G$ be any \infdef commutative divisible
  group in $L$. Then $G$ is definably isomorphic to some $\Apinf :=
 \bigcap_n p^n A(L)$, for $A$ a semi-Abelian variety defined over
  $L$.

2. If $A$ is a semi-Abelian variety defined over $L$, $\Apinf$,
which is the maximal  divisible subgroup of $A(L)$ is also the
smallest \infdef subgroup of $A(L)$ which is Zariski dense in
$A$.
\end{proposition}

Finally, this analysis, together with some results from \cite{delon2} and \cite{Hr1}, yields the full description
of minimal groups.

Before stating the actual  result, let us give some last definitions. The group $G$ is said to be {\em of linear
type} if for every $n$,  every definable subgroup of $G^{\times n}$ is a finite Boolean combination of translates
of definable subgroups of $G^{\times n}$. We define the {\it transcendence rank over $K$} of a group $G$, defined
over  $K$, to be the maximum of $\{\hbox {tr.degree}(K(g_\infty),K): g \in G\}$.

\vspace{-0.3cm}
\begin{proposition}\label{main} Let $G$ be an \infdef minimal group
  in $L$.
\begin{enumerate}
\item Either $G$ is not of linear type and then,
\begin{itemize}
\item  $G$ is definably isomorphic to the multiplicative
  group $({(\Lpinf)}^{*} , \cdot)$,
\item or $G$ is definably isomorphic to $E(\Lpinf)$ for $E$ an
  elliptic curve defined over $\Lpinf$,
\item or $G$ is definably isogenous to $(\Lpinf, +)$.
(isogenous here cannot be replaced by isomorphic).
\end{itemize}
\item Or  $G$ is of linear type and then,
\begin{itemize}
\item $G$ is divisible and $G$ is definably isomorphic to
  $\Apinf$ for some simple Abelian variety $A$ defined over $K$
  which is not isogenous  to an Abelian variety defined over
  $\Lpinf$,
\item or $G$ is of exponent $p$ and is definably isomorphic to a minimal
  \lclosed subgroup of $(L,+)$.
\end{itemize}
  In the divisible case $G$ has finite transcendence
  rank; in the exponent $p$ case, all transcendence
  ranks  are possible.
\end{enumerate}
\end{proposition}

The induced module-type structure on the minimal groups of exponent $p$ and of linear type is analyzed in
\cite{blos}.

A short word about some of the tools involved in the proofs of
Propositions \ref{defgroup} and \ref{main}: the proofs of \ref{defgroup},
\ref{diff} and \ref{gen} all involve at some point the classical
theorem of Weil's constructing an algebraic group from a generic
group law on a variety, or some generalizations of this theorem to
an abstract model theoretic context. In the specific case of separably closed
fields, another  fundamental tool   is the analysis of the
properties of the
{\em $\Lambda_n$-functors}, naturally associated to  the maps
$\lambda_n$: for each $n$, $\Lambda_n$ is a covariant functor from
the category of varieties $V$ defined over $K$ to itself,  with the
property  that the $L$-rational points of the variety
$\Lambda_n V$ are exactly the image by the map $\lambda_n$ of the
$L$-rational points of $V$.  In the case
of an algebraic group defined over $K$,
$\Lambda_1$ is equal to the
composition of the inverse of
the
Frobenius and of the classical Weil restriction of scalars
functor from $K^{1/p}$ to $K$.

Finally, the way we have stated Proposition \ref{main} uses the fact that
if a minimal group is not of linear type, then it is non orthogonal
to $\Lpinf$ (and hence definably isogenous to the
$\Lpinf$-rational points of some definable group over
$\Lpinf$). The only known proof of this so far uses the
powerful abstract machinery of  Zariski structures from
\cite{HrZi}.
This dichotomy result, for the particular case of groups of the
form $\Apinf$, is essential in  Hrushovski's proof of the Mordell-Lang
conjecture in characteristic $p$, which is still the only existing proof for
the general case.
In a recent paper Pillay and Ziegler
(\cite{PiZi}),
show that, with some extra
assumptions on
$A$, one can replace in this proof the heavy Zariski structure argument
by a much more elementary one.
These extra assumptions are satisfied when $A$ is  an {\em
  ordinary} semi-Abelian
variety (i.e. $A$ has the maximum possible number of $p^n$-torsion
points for every $n$), case which  was already
covered by previous non model-theoretic proofs (see \cite{AbVo}).

\subsection{Final remarks and questions}\label{section4}

\vskip-5mm \hspace{5mm}

As we have already mentioned earlier, the groups of finite
dimension definable in these ``enriched''  theories of fields  play
a major role in the applications of model theory to Diophantine
geometry. In the characteristic zero case,  the relevant groups are
the definable
subgroups of the group of rational points of Abelian varieties  in
differentially closed fields (Mordell-Lang conjecture for function
fields \cite{Hr1}), in generic difference fields (the Manin-Mumford
conjecture \cite{Hr2}, \cite{B2} and
the Tate-Voloch conjecture for semi-Abelian varieties defined over
${\mathbb Q}_p$ \cite{scanlon1}, \cite{scanlon2}). In
the characteristic $p$ case, the relevant groups are: the
$\infty$-definable divisible subgroups of the group of rational points of
semi-Abelian varieties
in separably closed
fields (the Mordell-Lang conjecture for function fields \cite{Hr1})
and
the definable subgroups of  the additive groups in
generic difference fields of characteristic $p$  (Drinfeld modules
\cite{scanlon3}).

One should note that,  in fact, separably closed fields are just
another instance of a {\em field with extra operators} (derivations or
automorphisms): one can equip any
separably closed field $L$ of finite degree of imperfection,
with an infinite family of Hasse derivations in such a way that the
resulting structure is bi-definably equivalent with $L$ considered as a
structure in the language described in section \ref{language}.
There are many interesting other possible types of ``enriched''
fields in this sense where the complete analysis of the model theoretic
structure remains to be done.

Finally,  one crucial step towards possible further
applications of the fine
study of finite rank definable sets to geometry would be an
understanding of the structure induced on the
so-called {\em trivial} or {\em disintegrated}
definable (or infinitely
definable) minimal  sets, that is the minimal sets such that the
induced pregeometry is disintegrated. This condition immediately
rules out definable groups. The absence of any well-understood
algebraic structure living on these ``trivial''  sets makes them
very difficult
to analyze.  The only results obtained so
far are in the context of differentially closed fields of
characteristic $0$: Hrushovski (\cite{udijoan}), building on some
results of Jouanolou (\cite{joan}), showed that in any trivial strongly
minimal set defined by a differential equation of order one,
the induced pregeometry is locally finite. The question of whether
this is true for higher order equations is still open.

\end{document}